\documentclass{amsart}

\usepackage{amsthm,amsmath,amsthm,amssymb,mathrsfs,graphicx,mathtools,relsize,hyperref,enumerate}
\usepackage[dvipsnames]{xcolor}
\usepackage[utf8]{inputenc}
\usepackage[all]{xy}

\newtheoremstyle{myplain}
{11pt}
{0pt}
{}
{0pt}
{\bfseries}
{.}
{ }
{\thmname{#1}\thmnumber{ #2}\thmnote{\textnormal{{ (#3) }}}}

\newtheoremstyle{myfancy}
{11pt}
{0pt}
{\itshape}
{0pt}
{\bfseries}
{.}
{ }
{\thmname{#1}\thmnumber{ #2}\thmnote{\textnormal{{ (#3) }}}}

\theoremstyle{myplain}

\newtheorem{defn}{Definition}[section]
\newtheorem{ex}[defn]{Example}
\newtheorem*{preuve}{Proof}
\newtheorem*{rmq}{Remark}

\theoremstyle{myfancy}

\newtheorem{prop}[defn]{Proposition}
\newtheorem{theo}[defn]{Theorem}
\newtheorem*{introthm}{Theorem}

\newtheorem{lemme}[defn]{Lemma}

\newcommand{\set}[1]{\left\{ #1 \right\}}

\newcommand{\Z}{\mathbb{Z}}

\newcommand{\Q}{\mathbb{Q}}
\newcommand{\F}{\mathbb{F}}

\title[Infinite families of triangle presentations]{Infinite families of triangle presentations}

\author[Alex Loué]{Alex Loué} 
\address{Institut de recherche en mathématique et physique \\
Chemin du Cyclotron 2 \\
boîte L7.01.02 \\
Université catholique de Louvain \\ 
1348 Louvain-la-Neuve \\
Belgique.}
\email{alex.loue@uclouvain.be}

\begin{document}

\begin{abstract}

A triangle presentation is a combinatorial datum that encodes the action of a group on a $2$-dimensional triangle complex with prescribed links, which is simply transitive on the vertices. We provide the first infinite family of triangle presentations that give rise to lattices in exotic buildings of type $\widetilde{\text{A}_2}$ of arbitrarily large order. Our method also gives rise to infinite families of triangle presentations for other link types, such as opposition complexes in Desarguesian projective planes.

\end{abstract}

\maketitle
\tableofcontents


\section{Introduction}

The construction and classification of groups acting simply transitively on the vertices of triangle buildings has seen two major contributions. First in \cite{cartwright1} and \cite{cartwright2}, this has been done for buildings of type $\widetilde{\text{A}_2}$ of orders $2$ and $3$. Second in \cite{carbone}, this has been done for triangle buildings such that the link of every vertex is isomorphic to the smallest finite thick generalized quadrangle, which has order $2$. In both cases, the approach consists in enumerating all triangle presentations compatible with the given datum, which amounts to constructing finite cell complexes whose universal cover is the desired building. There are various sources of motivation for constructing and studying such groups. Let us mention a few of interest to us.

In type $\widetilde{\text{A}_2}$, such groups are subject to the Normal Subgroup Theorem of \cite{bader1}, which states that every non-trivial normal subgroup must have finite index. In particular, any such group obeys the following strong dichotomy : either it is residually finite, or it is virtually simple. Of course, arithmetic lattices in Bruhat-Tits buildings are linear, and thus residually finite. In contrast, Galois lattices (i.e. lattices in Bruhat-Tits buildings that are not virtually contained in the linear part of the automorphism group) and lattices in non-Bruhat-Tits buildings are shown in \cite{bader2} to be non-linear. The latter buildings and lattices are called exotic. We do not know if Galois lattices exist. It is conjectured that exotic lattices are virtually simple.

Besides the aforementioned examples that result from the enumeration in \cite{cartwright1} and \cite{cartwright2}, the same authors describe an infinite family of such groups, one for each prime power order $q\geq 2$. To our knowledge, this is currently the only known infinite family of groups acting simply transitively on the vertices of buildings of type $\widetilde{\text{A}_2}$. Each group in this family embeds as an arithmetic lattice in $\text{PGL}_3(\F_q((t)))$.

Another construction of infinitely many lattices on $\widetilde{\text{A}_2}$ buildings was carried out in \cite{essert}, see also \cite{witzel} for a study of their (non-) arithmeticity. These lattices are regular on panels, and it turns out that an overwhelming majority of them are exotic. We also mention \cite{radu}, where a triangle presentation is given for a vertex regular building whose link is non-Desarguesian ; indeed, the link is the Hughes plane of order $9$. Such a building, and any lattice in it, must be exotic.

In this paper, we investigate triangle presentations under a more general perspective, in view of constructing new families of interesting examples of groups acting simply transitively on non-positively curved triangle complexes, possibly with other link types. More specifically, we study triangle presentations which display a large amount of symmetry. In this way, we are able to provide a construction of groups acting simply transitively on the vertices of residually Desarguesian, exotic buildings of type $\widetilde{\text{A}_2}$ of arbitrarily large order. To our knowledge, this is the first such infinite family of examples. We emphasize the following result :

\begin{introthm}

Let $p$ be a prime and let $q=p^e$, $e\geq 1$. Then, up to isomorphism, there exist at least : \[ \frac{2^{R(q)}-e(q-1)q(q+1)}{2e(q-1)^2 q^3 (q+1) } , \] quasi-isometry classes of groups acting simply transitively on the vertices of residually Desarguesian exotic buildings of type $\widetilde{\emph{A}_2}$ and order $q$, where :  \[  R(q) = \begin{cases} \frac{q+1}{3} &  \emph{if $q=-1$ mod $3$}, \\ \frac{q}{3} &  \emph{if $q=0$ mod $3$}, \\ \frac{q-1}{3} &  \emph{if $q=1$ mod $3$}.  \end{cases} \]

\end{introthm}

On the other hand, the groups constructed in \cite{carbone} are hyperbolic, and for those also it is interesting to investigate residual finiteness. By the recent work of \cite{ashcroft}, it has been shown that those groups are in fact virtually special, i.e. they act properly on non-positively curved, special cube complexes, for which the famous linearity theorem of \cite{agol} applies. A known obstruction to the existence of a proper action on a non-positively curved cube complex is Kazhdan's property (T), and this motivates the study of hyperbolic groups with Kazhdan's property (T) that arise from triangle presentations for other link types. In the same spirit as before, we are able to prove the following :

\begin{introthm}

Let $q\geq 2$ be a prime power. If $q=1\emph{ mod } 3$, then there exist at least $2^{\frac{q-1}{3}}$ distinct triangle presentations that give rise to groups acting simply transitively on the vertices of a $2$-dimensional, $\emph{CAT}(0)$ triangle complex such that the link of every vertex is isomorphic to the opposition complex of the Desarguesian projective plane of order $q$.

\end{introthm}

For $q\geq 5$, the groups that arise in this way have Kazhdan's property (T). We do not know if such groups are hyperbolic, but we believe that some of them could be. As mentioned above, such groups would not be virtually special. Either way, they are not known to be residually finite, nor linear. Their geometry is closely related to that of the Kac-Moody-Steinberg groups of type $\widetilde{\text{A}_2}$ introduced in \cite{caprace}.

In Section \ref{sec1}, we give a general definition of triangle presentations and an appropriate notion of equivalence between them, as well as further properties of groups that arise from them. In Section \ref{sec2}, we provide a construction of symmetric triangle presentations under the given of a finite group $G$ and a subset $S\subseteq G$ in special configuration. In Section \ref{sec3}, we apply this construction to obtain symmetric triangle presentations supported on the incidence graph of Desarguesian projective planes, and study the (non-) arithmeticity of the lattices that arise from them. In Section \ref{sec4} we provide some further examples of symmetric triangle presentations, among which the class of groups that we have just mentioned in the previous paragraph.

The author would like to acknowledge the work of Amy Herron, whose contributed talk at the \emph{Buildings $2023$} conference was a source of inspiration for studying unusually symmetric triangle presentations. We recommend interested readers to look out for the upcoming preprint \cite{herron}. The author would also like to thank Pierre-Emmanuel Caprace for helpful comments on an earlier version of this paper and support all troughout this project.


\section{Generalities on triangle presentations}
\label{sec1}

In this section, we define a general notion of triangle presentation, which gives rises to the simply transitive action of a group on the vertices of a triangle complexe with the property that the link of every vertex is prescribed.

\begin{figure}
\begin{center}
\includegraphics[scale=0.7]{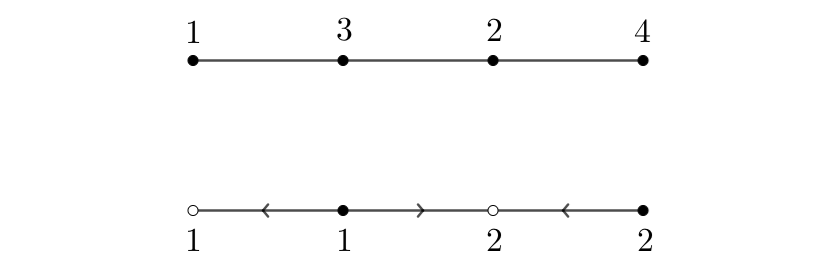}
\caption{The graph $\mathcal{G}^{(2)}_F$ constructed from $F= \set{(1,1),(2,1),(2,2)}$. Formally, it is the first one, but as explained we will use the second kind of representation all throughout this paper.}
\label{line515}
\end{center}
\end{figure}

Let $n\geq 1$ be an integer. For every $F\subseteq \set{1,\dots,n}^2$, we denote by $\mathcal{G}^{(n)}_F$ the simple graph with vertex set $\set{1,\dots,2n}$ and edge set : \[ \set{ \set{i,j+n} : (i,j) \in F } . \] It affords the obvious bipartition : \[ \set{1,\dots,2n} = \set{1,\dots,n} \cup \set{n+1,\dots,2n} . \] We shall think of the edges as being oriented from the largest to the smallest vertex ; the vertices can then be numbered with the same index set $\set{1,\dots,n}$ by identifying $k\in \set{n+1,\dots,2n}$ with $k-n$. See Figure \ref{line515} for a simple example.

Note that the map : \[  i \in \set{1,\dots,n} \mapsto i+n \in \set{n+1,\dots,2n} , \] defines a bijection between the vertices in each class of the natural bipartition of $\mathcal{G}^{(n)}_F$. In \cite{cartwright2}, the choice of a bijection between points and lines of a finite projective plane is the crucial parameter for the definition of a triangle presentation. In our approach, this bijection is the same for all representations of a given graph of the form $\mathcal{G}^{(n)}_F$, and the defining parameter becomes the set $F$ itself. Of course, both points of view are equivalent, but we have found it easier to work in this setting.

\begin{defn}

Let $F\subseteq \set{1,\dots,n}^2$. We say that $T\subseteq \set{1,\dots,n}^3$ is a triangle presentation compatible with $F$ if the following properties are satisfied : \begin{enumerate} \item If $(i,j,k) \in T$, then $(i,j) \in F$. \item If $(i,j) \in F$, then there exists a unique $k\in \set{1,\dots,n}$ such that $(i,j,k) \in T$. \item If $(i,j,k) \in T$, then $(j,k,i) \in T$. \end{enumerate}

\end{defn}

\begin{rmq}

For a given $F\subseteq \set{1,\dots,n}^2$, there does not necessarily exist compatible triangle presentations. See Figure \ref{line515} for such an example.

\end{rmq}

If $T\subseteq \set{1,\dots,n}^3$ is a triangle presentation compatible with a given $F\subseteq \set{1,\dots,n}^2$, we shall be interested in the group defined by the following presentation : \[ \Gamma^{(n)}_T = \langle a_l : l\in \set{1,\dots,n} \vert a_ia_ja_k : (i,j,k) \in T \rangle . \] Its presentation complex $\Sigma^{(n)}_T$ consists of triangles whose vertices are all identified and whose edges are glued following the labelling provided by $T$. It has a natural metric that restricts to the euclidean metric on each triangle and assigns length $1$ to each edge. We denote its universal cover by $\Delta^{(n)}_T$ ; recall that the group $\Gamma^{(n)}_T$ acts on it by deck transformations, which is an isometric action for the pullback of the metric on $\Sigma^{(n)}_T$. Its $1$-skeleton identifies with the Cayley graph of $\Gamma^{(n)}_T$ with respect to the canonical generating set $A= \set{a_l : l\in \set{1,\dots,n}}$, on which the group acts regularly.

\begin{figure}

\begin{center}
\includegraphics[scale=0.5]{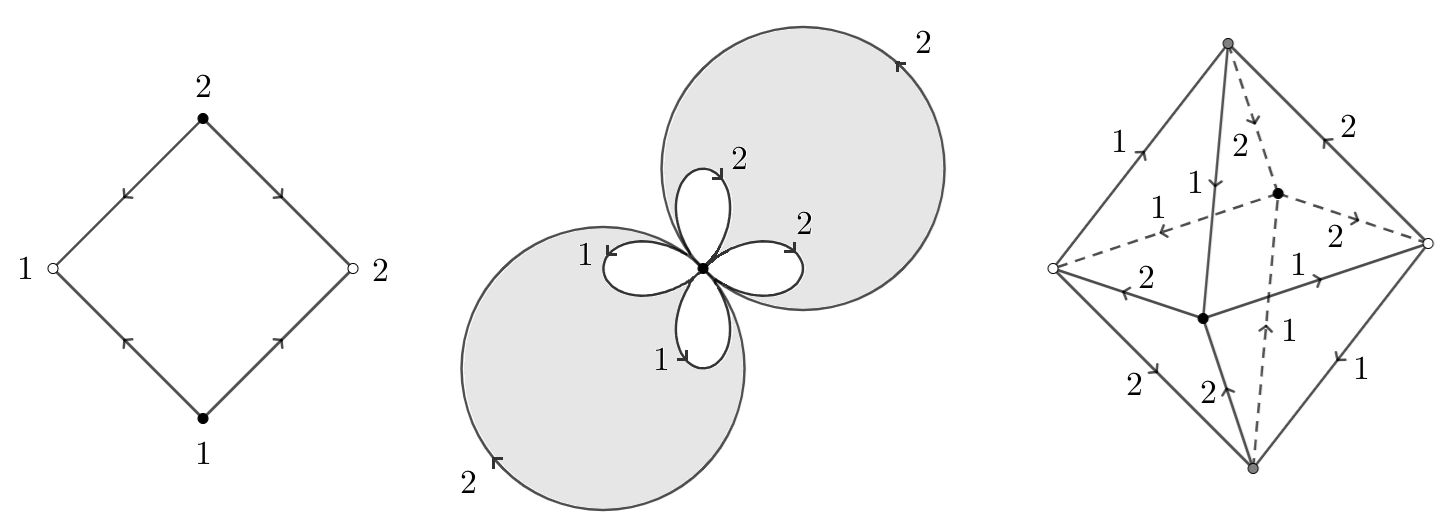}
\caption{From left to right, we have illustrated $\mathcal{G}^{(2)}_F$, $\Sigma^{(2)}_T$ and $\Delta^{(2)}_T$ from Example \ref{exoctahedron}, respectively.}
\label{octahedron1455}
\end{center}
\end{figure}

\begin{ex}
\label{exoctahedron}

Let : \[ F = \set{ (1,1),(1,2),(2,1),(2,2) } . \] Then : \[ T = \set{(1,1,2),(1,2,1),(2,1,1),(2,2,2)} , \] is a triangle presentation compatible with $F$. In Figure \ref{octahedron1455}, we have illustrated the graph $\mathcal{G}^{(2)}_F$, the cell complex $\Sigma^{(2)}_T$ and its universal cover $\Delta^{(2)}_T$, which has the structure of an octahedron. We make the observation that the link of every vertex in the octahedron is a square, which is isomorphic to the graph $\mathcal{G}^{(2)}_F$. It turns out that the group $\Gamma^{(2)}_T$ is cyclic of order $6$ in this case.

\end{ex}

\begin{theo}
\label{thmgeneral}

Let $F\subseteq \set{1,\dots,n}^2$. If $\mathcal{G}^{(n)}_{F}$ is connected and has girth at least $6$, and if $T\subseteq \set{1,\dots,n}^3$ is a triangle presentation compatible with $F$, then : \begin{enumerate} \item $\Delta^{(n)}_T$ is a $\emph{CAT}(0)$ triangle complex such that the link of every vertex is isomorphic to $\mathcal{G}^{(n)}_F$. \item There exists a $3$-coloring of the vertices of $\Delta^{(n)}_T$ such that the action of $\Gamma^{(n)}_T$ is type-rotating. In particular, $\Gamma^{(n)}_T$ has a normal subgroup of index $3$ which is torsion free. \item If $\mathcal{G}^{(n)}_F$ has girth at least $8$, then $\Gamma^{(n)}_T$ is hyperbolic. \item If the spectral gap of $\mathcal{G}^{(n)}_F$ is strictly greater than $\frac{1}{2}$, then $\Gamma^{(n)}_T$ has Kazhdan's property \emph{(T)}. \item If $\mathcal{G}^{(n)}_F$ is a generalized $m$-gon with $m\geq 3$, then $\Delta^{(n)}_T$ is a triangle building of type $(m,m,m)$. \end{enumerate}

\end{theo}

\begin{preuve}

Properties $1$ and $3$ follow from applications of the Cartan-Hadamard in the case of a $2$-dimensional polyhedral cell complex, see \cite{bridson} Part II, Chapter $5$. Property $2$ follows from considering the partition of the elements of $\Gamma^{(n)}_T$ by cosets of the kernel of the homomorphism $\Gamma^{(n)}_T \to \Z/3\Z$ defined by $a_l \mapsto 1$ for all $l\in \set{1,\dots,n}$. The kernel is the required normal subgroup of index $3$. It is torsion free because now the action is type-preserving, and we may apply Cartan's fixed point theorem, see \cite{bridson} Part II, Chapter $2$. Property $4$ is a direct application of Zuk's criterion, see \cite{Kazhdan} Section $5.6$. Property $5$ follows from \cite{tits}. $\square$

\end{preuve}

The remainder of this section is concerned with the classification of triangle presentations, up to a natural notion of equivalence. In the case of triangle presentations supported on generalized $3$-gons, this notion of equivalence actually accounts for all possible isomorphisms between associated groups.

Consider the wreath product $\text{Sym}(n) \wr \Z/2\Z$. We write the elements of this group in the form : \[ w = (\sigma_1,\sigma_2 ; e) , \] for all $\sigma_1,\sigma_2\in \text{Sym}(n)$ and $e\in \set{0,1}$. It acts on $\set{1,\dots,n}^2$ in the following way : the normal subgroup $\text{Sym}(n)\times \text{Sym}(n)$ acts component-wise, and the $\Z/2\Z$ factor acts via the reflexion $\rho$ defined by : \[ \rho(i,j) = (j,i) , \] for all $i,j\in \set{1,\dots,n}$. In turn, this action naturally induces an action on the set of subsets of $\set{1,\dots,n}^2$.

\begin{prop}

Let $F,F' \subseteq \set{1,\dots,n}^2$. Then, the graphs $\mathcal{G}^{(n)}_{F}$ and $\mathcal{G}^{(n)}_{F'}$ are isomorphic if and only if there exists $w\in \emph{Sym}(n)\wr \Z/2\Z$ such that $F' = wF$.

\end{prop}

\begin{preuve}

A permutation $\varphi : \set{1,\dots,2n} \to \set{1,\dots,2n}$ induces an isomorphism from $\mathcal{G}^{(n)}_{F}$ to $\mathcal{G}^{(n)}_{F'}$ if and only if $\varphi$ preserves the bipartition : \[ \set{1,\dots,2n} = \set{1,\dots,n} \cup \set{n+1,\dots,2n}, \] and if $\set{\varphi(i),\varphi(j+n)}$ is an edge of $\mathcal{G}^{(n)}_{F'}$ for all $(i,j) \in F$. There are two cases. If $\varphi(\set{1,\dots,n}) = \set{1,\dots,n}$, then if we let $\sigma_1$ be the restriction of $\varphi$ to $\set{1,\dots,n}$ and $\sigma_2(j) = \varphi(j+n)-n$, then $F'=(\sigma_1,\sigma_2 ; 0) F$. Else if $\varphi(\set{1,\dots,n}) = \set{n+1,\dots,2n}$, then if we let $\sigma_1(i) = \varphi(i)-n$ and $\sigma_2(j) = \varphi(j+n)$, then $F' = (\sigma_1,\sigma_2;1)F$. $\square$

\end{preuve}

Throughout the remainder of this section, let $F\subseteq \set{1,\dots,n}^2$.

\begin{prop}
\label{propdiag}

Let $\sigma \in \emph{Sym}(n)$ and let : \[ F' = \sigma F = \set{(\sigma(i),\sigma(j)) : (i,j) \in F} . \] If $T\subseteq \set{1,\dots,n}^3$ is a triangle presentation compatible with $F$, then : \[ T'  = \sigma T = \set{(\sigma(i),\sigma(j),\sigma(k)) : (i,j,k) \in T} ,  \] is a triangle presentation compatible with $F'$. Moreover, the assignments : \[ a_l \mapsto a'_{\sigma(l)} , \] where $l\in \set{1,\dots,n}$, extend to a group isomorphism $\Gamma^{(n)}_T \cong \Gamma^{(n)}_{T'}$.

\end{prop}

\begin{preuve}

Both statements are a direct verification of the defining properties of a triangle presentation, and that the relations are satisfied by the assignments. $\square$

\end{preuve}

\begin{prop}
\label{proprho}

Let : \[ F' = \rho F = \set{(j,i) : (i,j) \in F} . \] If $T\subseteq \set{1,\dots,n}^3$ is a triangle presentation compatible with $F$, then : \[ T' = \rho T = \set{ (j,i,k) : (i,j,k) \in T } , \] is a triangle presentation compatible with $F'$. Moreover, the assignments : \[ a_l \mapsto (a'_{l})^{-1} , \] where $l\in \set{1,\dots,n}$, extend to a group isomorphism $\Gamma^{(n)}_T \cong \Gamma^{(n)}_{T'}$.

\end{prop}

\begin{preuve}

Both statements are a direct verification of the defining properties of a triangle presentation, and that the relations are satisfied by the assignments. $\square$

\end{preuve}

\begin{defn}

Let $F'\subseteq \set{1,\dots,n}^2$. We say that $F'$ is equivalent to $F$ in a type-preserving way if $F' = \sigma F$ for some $\sigma \in \text{Sym}(n)$, and we say that $F'$ is equivalent to $F$ if it is equivalent to $F$ or $\rho F$ in a type-preserving way. If $T,T'\subseteq \set{1,\dots,n}^3$ are triangle presentations compatible with $F$ and $F'$ respectively, then we say that $T'$ is isomorphic to $T$ if $T' = \sigma T$ or $T' = \sigma \rho T$ for some $\sigma \in \text{Sym}(n)$.

\end{defn}

By Propositions \ref{propdiag} and \ref{proprho}, if $F'\subseteq \set{1,\dots,n}^2$ is equivalent to $F$, then all triangle presentations compatible with $F$ carry over to triangle presentations compatible with $F'$ and yield isomorphic groups. The stabiliser of $F$ in $\text{Sym}(n)$ is denoted by $\text{Aut}^+(F)$, and the stabiliser of $F$ in $\text{Sym}(n) \times \Z/2\Z$ (where the $\Z/2\Z$ factor acts via $\rho$) is denoted by $\text{Aut}(F)$. Moreover, the group $\text{Aut}(F)$ acts on the set of triangle presentations compatible with $F$. If $T\subseteq \set{1,\dots,n}^3$ is a triangle presentation compatible with $F$, then the stabiliser of $T$ in $\text{Aut}(F)$ is denoted by $\text{Aut}(T)$, and we similarly define $\text{Aut}^+(T) = \text{Aut}(T) \cap \text{Aut}^+(F)$.

\begin{prop}
\label{propcount}

Let $F\subseteq \set{1,\dots,n}^2$ and let $R$ be a set of representatives of all isomorphism classes of triangle presentations compatible with $F$. Then, the total number of triangle presentations compatible with $F$ is equal to $ \sum_{T\in R}\frac{|\emph{Aut}(F)|}{|\emph{Aut}(T)|}$.

\end{prop}

\begin{preuve}

Apply the orbit-stabiliser theorem to the action of $\text{Aut}(F)$ on the set of triangle presentations compatible with $F$. $\square$

\end{preuve}

\begin{ex}
\label{exalt}

Let : \[ F = \set{(1,2),(1,3),(1,4),(2,1),(2,3),(2,4),(3,1),(3,2),(3,4),(4,1),(4,2),(4,3)}  . \] Clearly, since $(i,j) \in F$ if and only if $i\neq j$, we see that $\text{Aut}^+(F) = \text{Sym}(4)$ can be extended by $\rho$, and thus $\text{Aut}(F) = \text{Sym}(4) \times \Z/2\Z$. An exhaustive enumeration shows that there are exactly two triangle presentations compatible with $F$, namely : \[ T_1 =  \set{(1,2,3),(1,3,4),(1,4,2),(2,1,4),(2,3,1),(2,4,3),(3,1,2),(3,2,4),(3,4,1),(4,1,3),(4,2,1),(4,3,2)}  , \] \[ T_2 = \set{(1,2,4),(1,3,2),(1,4,3),(2,1,3),(2,3,4),(2,4,1),(3,1,4),(3,2,1),(3,4,2),(4,1,2),(4,2,3),(4,3,1)}   . \] We have $T_2 = (1)(2)(3,4) T_1$, so they are isomorphic and $\text{Aut}(T_1)$ must have index $2$ in $\text{Aut}(F)$. In fact, one can check that $\text{Aut}^+(T_1) = \text{Alt}(4)$ can be extended by $(1)(2)(3,4)\rho$, and therefore $\text{Aut}(T) = \text{Alt}(4) \rtimes \Z/2\Z$ where the $\Z/2\Z$ factor acts on $\text{Alt}(4)$ as conjugation by $(1)(2)(3,4)$.

\end{ex}

We can exploit the quasi-isometry rigidity of buildings of type $\widetilde{\text{A}_2}$ to show the following result.

\begin{prop}
\label{propqi}

Let $F_1,F_2\subseteq \set{1,\dots,n}^2$. If both $\mathcal{G}^{(n)}_{F_i}$, $i=1,2$ are isomorphic to the same generalized $3$-gon and if $T_i \subseteq \set{1,\dots,n}^3$ are triangle presentations compatible with $F_i$, then the $\Gamma^{(n)}_{T_i}$ are quasi-isometric if and only if the $T_i$ are isomorphic.
    
\end{prop}

\begin{preuve}

If $T_1$ is isomorphic to $T_2$, then we have seen in Proposition \ref{propdiag} and Proposition \ref{proprho} that there is a natural isomorphism between the $\Gamma^{(n)}_{T_i}$, so we focus on the converse. If the $\Gamma^{(n)}_{T_i}$ are quasi-isometric, then the associated buildings $\Delta^{(n)}_{T_i}$ are quasi-isometric and by the main result of \cite{kleinerleeb}, this implies that the buildings are isometric. In particular, since the lattices $\Gamma^{(n)}_{T_i}$ are vertex transitive, we may assume without loss of generality that there is an isometry $\varphi : \Delta^{(n)}_{T_1} \to \Delta^{(n)}_{T_2}$ such that $\varphi$ carries the identity element in $\Gamma^{(n)}_{T_1}$ to the identity element in $\Gamma^{(n)}_{T_2}$. In particular, $\varphi$ induces a permutation of the generators which defines an isometry of the links, and thus an isomorphism between $T_1$ and $T_2$. $\square$
    
\end{preuve}


\section{Symmetric triangle presentations}
\label{sec2}

In this section, we will use the fact that some graphs afford a regular action of a subgroup of their automorphism group, which induces an intrinsic numbering of the vertices. This leads to an explicit criterion for the existence of symmetric triangle presentations.

Let $n\geq 1$ and let $G=\set{g_1,\dots,g_n}$ be a finite group of order $n$. For every $S\subseteq G$, let : \[ F(G,S) = \set{(x,xs) : x\in G,s\in S} . \] When we write $\mathcal{G}^{(n)}_{F(G,S)}$, we mean $\mathcal{G}^{(n)}_F$ where : \[ F = \set{ (i,j) : g_{i}^{-1}g_j \in S }  .  \] The equivalence class of $F$ constructed in this way does not depend on the numbering of the elements of $G$, and so we will allow ourselves this abuse of notation. Similarly, if $T\subseteq \set{1,\dots,n}^3$ is a triangle presentation compatible with $F$, we will identify it with : \[ T = \set{ (g_i,g_j,g_k) : (i,j,k) \in T } , \] and say that it is a triangle presentation compatible with $F(G,S)$. The diagonal action of $G$ by multiplication on the left preserves $F(G,S)$, and thus $G\leq \text{Aut}^+(F(G,S))$. Recall from Section \ref{sec1} the definition of the reflexion $\rho$, which in this setting takes on the form : \[ \rho(x,y) = (y,x), \] for all $x,y\in G$.

\begin{prop}
\label{propmu}

Let $\mu \in \emph{Sym}(G)$ be the permutation defined by : \[ \mu(x) = x^{-1} .  \] If $G$ is abelian, then $\mu \rho \in \emph{Aut}(F(G,S))$.

\end{prop}

\begin{preuve}

For all $x\in G$ and $s\in S$, we compute : \[ \mu \rho (x,xs) = \mu(xs,x) = (s^{-1}x^{-1},x^{-1}) =(x^{-1}s^{-1},x^{-1}s^{-1}s) = (y,yt) , \] where $y = x^{-1}s^{-1} \in G$ and $t = s\in S$. $\square$

\end{preuve}

\begin{prop}
\label{proplambda}

Let $S\subseteq G$. If $\lambda : S \to S$ is a map such that : \[ s\lambda(s)\lambda^2(s) = 1_G, \] for all $s\in S$, then : \[ T = \set{(x,xs,xs\lambda(s)) : x\in G , s\in S } , \] is a triangle presentation compatible with $F(G,S)$ such that $G\leq \emph{Aut}^+(T)$.

\end{prop}

\begin{preuve}

It suffices to verify the defining properties of a triangle presentation. The first two are obvious because $(x,y,z) \in T$ if and only if $z = y\lambda(x^{-1}y)$. The third follows from the following computation : \[ (xs,xs\lambda(s),x) = (xs,xs\lambda(s),xs\lambda(s)\lambda^2(s)) = (y,yt,yt\lambda(t)) , \] where $y = xs \in G$ and $t=\lambda(s)\in S$. $\square$ 

\end{preuve}

\begin{rmq}

In fact, it is quite easy to show the converse, namely that a $G$-invariant triangle presentation must be of this form, for some $\lambda : S \to S$. This is in fact how we were lead to the condition that $s\lambda(s)\lambda^2(s) = 1_G$.

\end{rmq}

\begin{ex}

Let $G = \Z/2\Z \times \Z/2\Z$ and let $S = \set{(1,0),(1,1),(0,1)}$. Let $\lambda : S \to S$ be defined via the assignments : \[ (1,0) \mapsto (0,1) , \text{ } (1,1) \mapsto (1,0) , \text{ } (0,1) \mapsto (1,1) . \] Note that it defines a permutation $\lambda : S \to S$ of order $3$. For all $s\in S$, we have indeed : \[ s + \lambda(s) + \lambda^2(s) = (0,0) , \] because all of the $s\in S$ belong to the same orbit of $\lambda$, and the sum of all non-trivial elements in $G$ equals the identity in this case. Therefore, by Proposition \ref{proplambda} we have that : \[ T = \set{(x,x+s,x+s+\lambda(s)) : x\in G , s\in S} , \] is a triangle presentation compatible with $F(G,S)$ and $G\leq \text{Aut}^+(F(G,S))$. In fact, $F(G,S)$ is equivalent to Example \ref{exalt}, and $T$ is isomorphic to the triangle presentations given there, and the subgroup $G\leq \text{Aut}^+(T)$ corresponds to the subgroup of $\text{Alt}(4)$ generated by products of disjoint transpositions.

\end{ex}

\begin{ex}

If $S\subseteq G$ is such that $s^3=1_G$ for all $s\in S$, then we can take $\lambda = \text{id}_S$ in Proposition \ref{proplambda}. The simplest such example (with $S$ non-empty) is for $G = \Z/3\Z$ and $S = \set{1,2}$, in which case : \[ F(G,S) = \set{(0,1),(0,2),(1,0),(1,2),(2,0),(2,1) } , \] and : \[ T = \set{(x,x+s,x+2s) : x\in G, s\in S} \] \[ = \set{(0,1,2),(0,2,1),(1,0,2),(1,2,0),(2,0,1),(2,1,0) } , \] is a compatible triangle presentation.

\end{ex}

Interestingly, under some further assumptions, the $G$-invariant triangle presentation of \ref{proplambda} affords some flexibility that allows for the contruction of many triangle presentations. This was a suprising discovery in our investigation of symmetric triangle presentations.

\begin{prop}
\label{propbreak}

Let $S\subseteq G$ and let $H\leq G$ be a subgroup. If $\lambda : S \to S$ is a permutation such that $\lambda^3=\emph{id}_S$ and : \[ s\lambda(s)\lambda^2(s) = s\lambda^2(s)\lambda(s) = 1_G, \] for all $s\in S$, let : \[ O = \set{ \emph{orbits of $\lambda$ of length $3$ which are contained in $H$} } , \] and let $o : S \mapsto \set{s,\lambda(s),\lambda^2(s)} \in \mathcal{P}(S)$. Then, for every map $\kappa : G/H \times O \to \set{\pm 1}$, we have that : \[ T_\kappa = \set{(x,xs,xs\lambda(s)) : x\in G , s\in S \emph{ with } o(s)\notin O } \] \[ \bigcup \set{ (x,xs,xs\lambda^{\kappa(xH,o(s))}(s)) : x\in G , s\in S \emph{ with } o(s) \in O } , \] is a triangle presentation compatible with $F(G,S)$. In particular, there are at least $2^{[G:H] |O|}$ distinct triangle presentations compatible with $F(G,S)$.

\end{prop}

\begin{rmq}

If $H\leq G$ is a normal subgroup, then one can check that $H\leq \text{Aut}^+(T_\kappa)$. For instance, this is the case for $H=G$.

\end{rmq}

\begin{preuve}

It suffices to verify the defining properties of a triangle presentation. As in the proof of Proposition \ref{proplambda}, the first two are obvious. The third is settled in two steps. Let $x\in G$ and $s\in S$. If $o(s) \notin O$, let $y = xs\in G$ and $t = \lambda(s)$ and note that $o(t)=o(s) \notin O$. Then : \[ (xs,xs\lambda(s),x) = (xs,xs\lambda(s),xs\lambda(s)\lambda^2(s)) = (y,yt,yt\lambda(t)) \in T_\kappa . \] On the other hand, if $o(s) \in O$, let $y = xs\in G$ and $t = \lambda^{\kappa(xH,o(s))}(s) \in S$. Observe that $o(s) \in O$ implies that $s\in H$, and thus $yH = xsH = xH$. Moreover, $o(t)=o(s)\in O$, and therefore : \[ \lambda^{\kappa(yH,o(t))}(t) = \lambda^{\kappa(xH,o(s))}(\lambda^{\kappa(xH,o(s))}(s)) = \lambda^{2 \kappa(xH,o(s))}(s) ,\] and thus : \[ (xs,xs\lambda^{\kappa(xH,o(s))}(s),x) = (xs,xs\lambda^{\kappa(xH,o(s))}(s),xs\lambda^{\kappa(xH,o(s))}(s)\lambda^{2 \kappa(xH,o(s))}(s)) \] \[ = (y,yt,yt\lambda^{\kappa(yH,o(t))}(t)) \in T_\kappa , \] which ends our proof. $\square$

\end{preuve}

\begin{ex}
\label{exquad}

Let $G = \Z/21\Z$ and $S = \set{7,9,14,15,18}$. The automomorpishm $x\in G \mapsto 4x \in G$ preserves the set $S$ ; if we denote its restriction by $\lambda : S \to S$, then one can check that $\lambda = (7)(9,15,18)(14)$, and in particular $\lambda^3 = \text{id}_S$ and : \[ s+\lambda(s)+\lambda^2(s) = s + \lambda^2(s) + \lambda(s) = (1+4+4^2)s = 21s = 0 , \] for all $s\in S$. Let $H\leq G$ be the subgroup of order $7$, which is generated by $3$. We note that $S\cap H = \set{9,15,18}$ consists of precisely one orbit of $\lambda$, and thus $O = \set{\set{9,15,18}}$. According to Proposition \ref{propbreak}, there are at least $2^{[G:H]|O|} = 2^3 = 8$ triangle presentation compatible with $F(G,S)$, and an exhaustive enumeration confirms that there are $8$ of them in total, hence they are all of the form $T_\kappa$ for some $\kappa : G/H\times O \to \set{\pm 1}$. In plain words, the triangle presentation $T_\kappa$ is obtained from the triangle presentation : \[ T = \set{(x,x+s,x+s+\lambda(s)) : x\in G , s\in S} , \] which we have seen in Proposition \ref{proplambda} is such that $G\leq \text{Aut}^+(T)$, and then inverting $\lambda$ along some orbits of length $3$ of $\lambda$, on some left-cosets of $H$ in $G$ ; which orbits are inverted and on which cosets is precisely the content of the map $\kappa$. Up to isomorphism, there are only two classes in this case, which are represtented by $T_1$ and $T_2$ given in Tables $1$ and $2$, and correspond to $\kappa$ equal to $1$ on all cosets of $H$ in $G$ and $\kappa$ equal to $1$ on all cosets of $H$ in $G$ except $2+H$, respectively. We emphasize that $T_1$ and $T_2$ share a good amount of common triples, and the differences in $T_2$ to $T_1$ are highlighted in red.

\end{ex}


\section{Construction of non-arithmetic $\widetilde{\text{A}_2}$-lattices}
\label{sec3}

In this section, we construct infinitely many triangle presenations which give rise to groups acting simply transitively on the vertices of buildings of type $\widetilde{\text{A}_2}$ for arbitrary prime power orders $q\geq 2$. After giving a simple arithmeticity criterion, we apply it to show that most of these lattices are exotic.

Let $p$ be a prime and let $q=p^e$ with $e\geq 1$. The simple graph which we denote by $\text{A}_2(\F_q)$ has vertex set the proper $\F_q$-subspaces of $(\F_q)^3$ and edge set : \[ \set{ \set{ V_1,V_2 } : \text{$V_1\subseteq V_2$ and $\text{dim}(V_i) = i$} } . \] It is a generalized $3$-gon of order $q$. The group $\text{P$\Gamma$L}_3(\F_q)$ acts faithfully and type preservingly by automorphisms on $\text{A}_2(\F_q)$, and has index $2$ in the full automorphism group, which has order $e(q^2+q+1)(q-1)^2q^3(q+1)$.

Let $G = \Z/(q^2+q+1)\Z$ and let : \[ U = \set{x\in \F_{q^3}^\times : \text{Tr}_{\F_q}(x) = 0 } . \]  Since the trace is linear, $U$ is preserved by $\F_q^\times$ and contains $q^2-1$ elements. We denote by $\overline{U}\subseteq \F_{q^3}^\times/\F_q^\times$ the image of $U$ in the quotient ; it contains $q+1$ elements. Let $\alpha \in \F_{q^3}$ be a primitive element. Then, there exists a unique isomorphism $\F_{q^3}^\times/\F_q^\times \cong G$ such that $\alpha$ identifies with $1$. Let $S\subseteq G$ be the subset that corresponds to $\overline{U}$ through that isomorphism. In other words : \[ S = \set{ l \in \set{0,\dots,q^2+q} : \text{Tr}_{\F_q}(\alpha^{l}) = 0 } . \]

\begin{prop}

The graph $\mathcal{G}^{(q^2+q+1)}_{F(G,S)}$ is isomorphic to $\emph{A}_2(\F_q)$.

\end{prop}

\begin{preuve}

We refer to \cite{singer} and \cite{pott} Section $2.1$, for a more recent account. $\square$

\end{preuve}

\begin{figure}
\begin{center}
\includegraphics[scale=0.5]{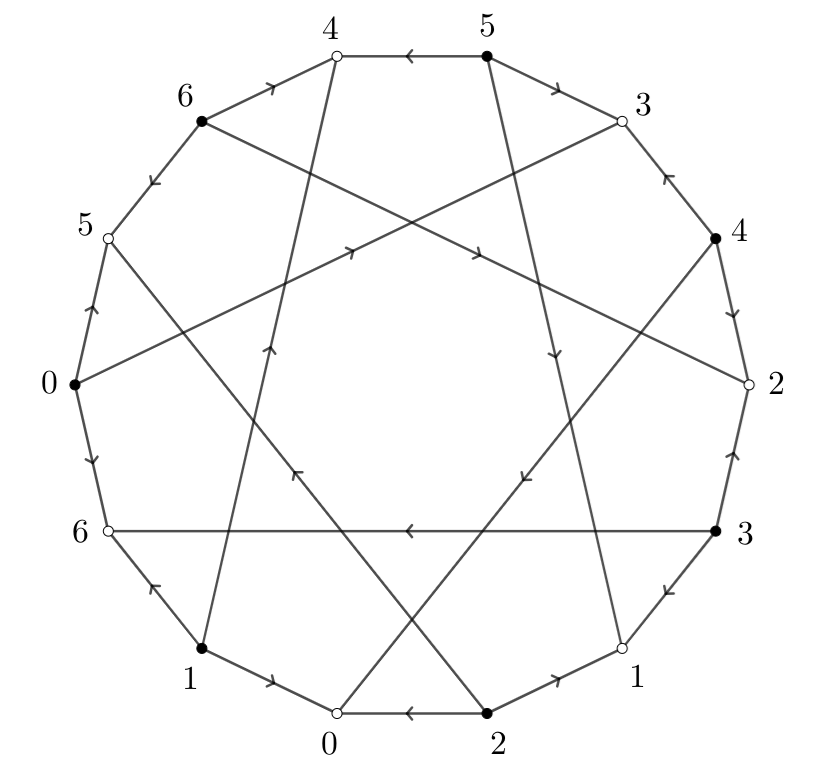}
\caption{A symmetric numbering of the Heawood graph, which is isomorphic to $\text{A}_2(\F_2)$.}
\label{heawood768}
\end{center}
\end{figure}

\begin{ex}
\label{exheawood}

For $q=2$, we have $G = \Z/7\Z$. There is a primitive element $\alpha \in \F_8$ with minimal polynomial $\alpha^3+\alpha + 1 = 0$ over $\F_2$ ; for that specific choice, we get $S = \set{1,2,4}$. We have illustrated the graph $\mathcal{G}^{(q^2+q+1)}_{F(G,S)}$ in Figure \ref{heawood768}.

\end{ex}

Let $\lambda : x\in G \mapsto qx \in G$. Since $q$ and $q^2+q+1$ are coprime, it defines an automorphism of $G$. It has order $3$, and more precisely, it satisfies : \[ x + \lambda(x) + \lambda^2(x) = (1+q+q^2)x = 0, \] for all $x\in G$.

\begin{prop}

The automorphism $\lambda : G \to G$ preserves $S$. In particular, $\lambda \in \emph{Aut}^+(F(G,S))$.

\end{prop}

\begin{preuve}

For all $x\in \F_{q^3}$, we have : \[ \text{Tr}_{\F_q}(x^q) = \text{Tr}_{\F_q}(x) . \] Moreover, the Frobenius automorphism $x\in \F_{q^3} \mapsto x^q \in \F_{q^3}$ fixes the subfield $\F_q$ pointwise, and thus it descends to an automorphism of $\F_{q^3}^\times/\F_q^\times$ which preserves $\overline{U}$. This statement carries over through the isomorphism $\F_{q^3}^\times/\F_q^\times \cong G$ and proves our claim. $\square$

\end{preuve}

 Let : \[ R(q) = \begin{cases} \frac{q+1}{3} &  \text{if $q=-1$ mod $3$}, \\ \frac{q}{3} &  \text{if $q=0$ mod $3$}, \\ \frac{q-1}{3} &  \text{if $q=1$ mod $3$}.  \end{cases} \]

\begin{prop}
\label{proporbit}

The restriction of $\lambda$ to $S$ has $R(q)$ orbits of length $3$.

\end{prop}

\begin{preuve}

Since $\lambda^3 = \text{id}_S$, it suffices to count the number of fixed points, or in other words, we need to count the number of solutions in $S$ to the equation : \[ (q-1)x = 0 \text{ mod } q^2+q+1 . \] Note that the gcd of $q-1$ and $q^2+q+1$ is equal to the gcd of $3$ and $q-1$, hence there are non-zero solutions in $G$ if and only if $q = 1$ mod $3$. Now we study each case separately. If $q=-1$ mod $3$, then $0$ is the only fixed point of $\lambda$, but $0 \notin S$ because : \[ \text{Tr}_{\F_q}(\alpha^0) = \text{Tr}_{\F_q}(1) = 3 \neq 0 .  \] Hence there are no fixed points in $S$, and thus there are $\frac{|S|}{3} = \frac{q+1}{3}$ orbits of length $3$ in $S$. If $q = 0$ mod $3$, then $0$ is the only fixed point of $\lambda$, and $0\in S$ because $\F_q$ has characteristic $3$, and thus : \[ \text{Tr}_{\F_q}(\alpha^0) = \text{Tr}_{\F_q}(1) = 3 = 0 . \] So there is one fixed point in $S$, and thus there are $\frac{|S|-1}{3} = \frac{q}{3}$ orbits of length $3$ in $S$. Finally, if $q=1$ mod $3$, then $0$ is a fixed point but $0\notin S$ for the same reason as before, and $\frac{q^2+q+1}{3}$ and $-\frac{q^2+q+1}{3}$ are the remaining non-zero fixed points. Moreover, we have : \[  \text{Tr}_{\F_q}(\alpha^{\frac{q^2+q+1}{3}}) = \alpha^{\frac{q^2+q+1}{3}} + (\alpha^{\frac{q^2+q+1}{3}})^q + (\alpha^{\frac{q^2+q+1}{3}})^{q^2} \] \[ = \alpha^{\frac{q^2+q+1}{3}}(1 + (\alpha^{\frac{q^2+q+1}{3}})^{q-1} + (\alpha^{\frac{q^2+q+1}{3}})^{q^2-1} ) \] \[ = \alpha^{\frac{q^2+q+1}{3}}( 1 + \beta + \beta^2 ) = 0 , \] where $\beta = \alpha^{\frac{q^3-1}{3}}$ satisfies $1+\beta+\beta^2 = 0$ because it has order $3$. This (and a similar computation for the opposite) shows that $\frac{q^2+q+1}{3}$ and $-\frac{q^2+q+1}{3}$ belong to $S$. So there are $2$ fixed points in $S$, and thus there are $\frac{|S|-2}{3} = \frac{q-1}{3}$ orbits of length $3$ in $S$. $\square$

\end{preuve}

Let $O$ be the set of orbits of $\lambda : S \to S$ which have length $3$ and let $o : s\in S \mapsto \set{s,\lambda(s),\lambda^2(s)} \in \mathcal{P}(S)$.

\begin{prop}

For every $\kappa : O \to \set{\pm 1}$, we have that : \[ T_\kappa = \set{(x,x+s,x+s+\lambda(s)) : x\in G, s\in S \emph{ with } o(s)\notin O } \] \[ \bigcup \set{ (x,x+s,x+s+\lambda^{\kappa(o(s))}(s)) : x\in G , s\in S \emph{ with } o(s) \in O } , \] is a triangle presentation compatible with $F(G,S)$ such that $G \leq \emph{Aut}^+(T_\kappa)$. In particular, there are at least $2^{R(q)}$ triangle presentations compatible with $F(G,S)$.

\end{prop}

\begin{preuve}

This is a direct application of Proposition \ref{propbreak}, with $H=G$. $\square$

\end{preuve}

\begin{rmq}
\label{remmu}

Recall from Proposition \ref{propmu} that $\mu \rho \in \text{Aut}(F(G,S))$. It is a direct computation to check that $T_{-\kappa} = \mu \rho T_{\kappa}$ for every $\kappa : O \to \set{\pm 1}$, so at most $2^{R(q)-1}$ of the triangle presentations constructed in this way are actually non-isomorphic.

\end{rmq}

\begin{ex}

For $q=2$, we have seen in Example \ref{exheawood} that $G=\Z/7\Z$ and $S=\set{1,2,4}$. We see that $S$ consists of one orbit of length $3$, and thus $O = \set{\set{1,2,4}}$, and there are two choices for $\kappa$. We provide the triangle presentation corresponding to $\kappa(o(1)) = 1$ in Table $3$. By Remark \ref{remmu}, the other triangle presentation is isomorphic to it. By exhaustive enumeration, one can check that there are no more triangle presentations compatible with $F(G,S)$.

\end{ex}

\begin{rmq}

When $\kappa : O \to \set{\pm 1}$ is the constant map equal to $1$, the triangle presentation $T_\kappa$ was originally found by Amy Herron, and further described in \cite{herron}, where it is shown that the corresponding lattice is arithmetic.
    
\end{rmq}

Let $\kappa : O \to \set{\pm 1}$, let $\Gamma = \Gamma^{(q^2+q+1)}_{T_\kappa}$ and let $\Delta = \Delta^{(q^2+q+1)}_{T_\kappa}$. The group $\Gamma$ acts regularly on the vertices of $\Delta$, and by Theorem \ref{thmgeneral}, we have that $\Delta$ is a triangle building of type $(3,3,3)$, to which we commonly refer to as a building of type $\widetilde{\text{A}_2}$. It is residually Desarguesian of order $q$, i.e. the link of every vertex is isomorphic to $\text{A}_2(\mathbb{F}_q)$. Let $v_0 = 1_\Gamma$ and $v_1 = a_0^{-1}$. They are adjacent in $\Delta$. Let $\Pi_i$ denote the link of $v_i$ in $\Delta$, in such a way that $v_{1-i} \in \Pi_i$. Let $\Lambda = \Pi_0 \cap \Pi_1$, and observe that : \[ \Lambda = \set{a_s : s\in S} . \] We further denote by $H_i$ the stabiliser of $v_{1-i}$ in $\text{Aut}(\Pi_i)$. It naturally preserves $\Lambda$, and we denote by $Q_i\leq \text{Sym}(\Lambda)$ the image of $H_i$ for the restriction of this action.

The remainder of this section builds up towards the proof of the main result of this paper.

\begin{theo}
\label{theomain}

Let $\sigma : a_s \in \Lambda \mapsto a_{\lambda^{\kappa(o(s))}(s)} \in \Lambda $. If $\sigma \notin Q_0$, then $\Delta$ is not isomorphic to a Bruhat-Tits building. Therefore, there are at least $2^{R(q)}-e(q-1)q(q+1)$ distinct choices for $\kappa$ such that $\Gamma$ is non-arithmetic, and those yield at least : \[ \frac{2^{R(q)} -e(q-1)q(q+1)}{2e(q-1)^2q^3(q+1)} , \] quasi-isometry classes of groups acting simply transitively on exotic buildings of type $\widetilde{\emph{A}_2}$ of order $q$.

\end{theo}

\begin{lemme}
\label{lemmaextension}

Let $v_2\in \Pi_0$ be opposite to $v_1$. Then, for all $h\in H_0$, there exists $h'\in H_0$ such that it has the same image as $h$ in $Q_0$ but $h'(v_2) = v_2$.

\end{lemme}

\begin{preuve}

Without loss of generality, we may assume that $\Pi_0$ is the simple graph $\text{A}_2(\F_q)$. Moreover, without loss of generality, by transitivity of its automorphism group, we may assume that : \[ v_1 = \langle \begin{bmatrix} 0 \\ 0 \\ 1  \end{bmatrix} \rangle , v_2 = \langle \begin{bmatrix} 1 \\ 0 \\ 0  \end{bmatrix} ,\begin{bmatrix} 0 \\ 1 \\ 0  \end{bmatrix}  \rangle . \] Then, we see that the Galois group of $\mathbb{F}_q$ over $\mathbb{F}_p$ fixes $v_1$ and $v_2$, hence it suffices to prove our claim for linear automorphisms. Thus, let : \[ h = \begin{bmatrix} A & 0 \\ * & *  \end{bmatrix} \in H_0 ,  \] where $A\in \text{GL}_2(\mathbb{F}_q)$. The image of $h$ in $Q_0$ is determined by the action of $A$ on the projective line. Hence, let : \[ h' = \begin{bmatrix} A & 0 \\ 0 & *  \end{bmatrix} \in H_0 , \] which has the same local action as $h$ but fixes $v_2$. $\square$

\end{preuve}

If $K$ is a field and $\nu$ is a discrete valuation on $K$, let $\mathcal{O}_K$ denote the valuation ring of $K$ with respect to $\nu$ and let $\pi \in K^\times$ be a uniformizer. If $e_1,e_2,e_3$ is a basis of $K^3$, we write : \[ L = \langle \langle e_1,e_2,e_3 \rangle \rangle , \] for the $\mathcal{O}_K$-lattice generated by the $e_i$. If $\pi L \subseteq L' \subseteq L$, we say that the lattices $L$ and $L'$ are incident. The simplicial complex which we denote by $\widetilde{\text{A}_2}(K,\nu)$ is the flag complex generated by the vertex set consisting of equivalence classes of $O_K$-lattices in $K^3$, with incidence described as above. It is a building of type $\widetilde{\text{A}_2}$, which we call the Bruhat-Tits building associated to $(K,\nu)$. The group $\text{GL}_3(K)$ acts on it by linear automorphisms. For more detail, see \cite{abramenko} Section $6.9$.

\begin{prop}
\label{proparithmetic}

Let $K$ be a field with discrete valuation $\nu$. If $\Delta$ is isomorphic to $\widetilde{\emph{A}_2}(K,\nu)$, then $Q_1 = Q_0$.

\end{prop}

\begin{figure}
\begin{center}
\includegraphics[scale=0.4]{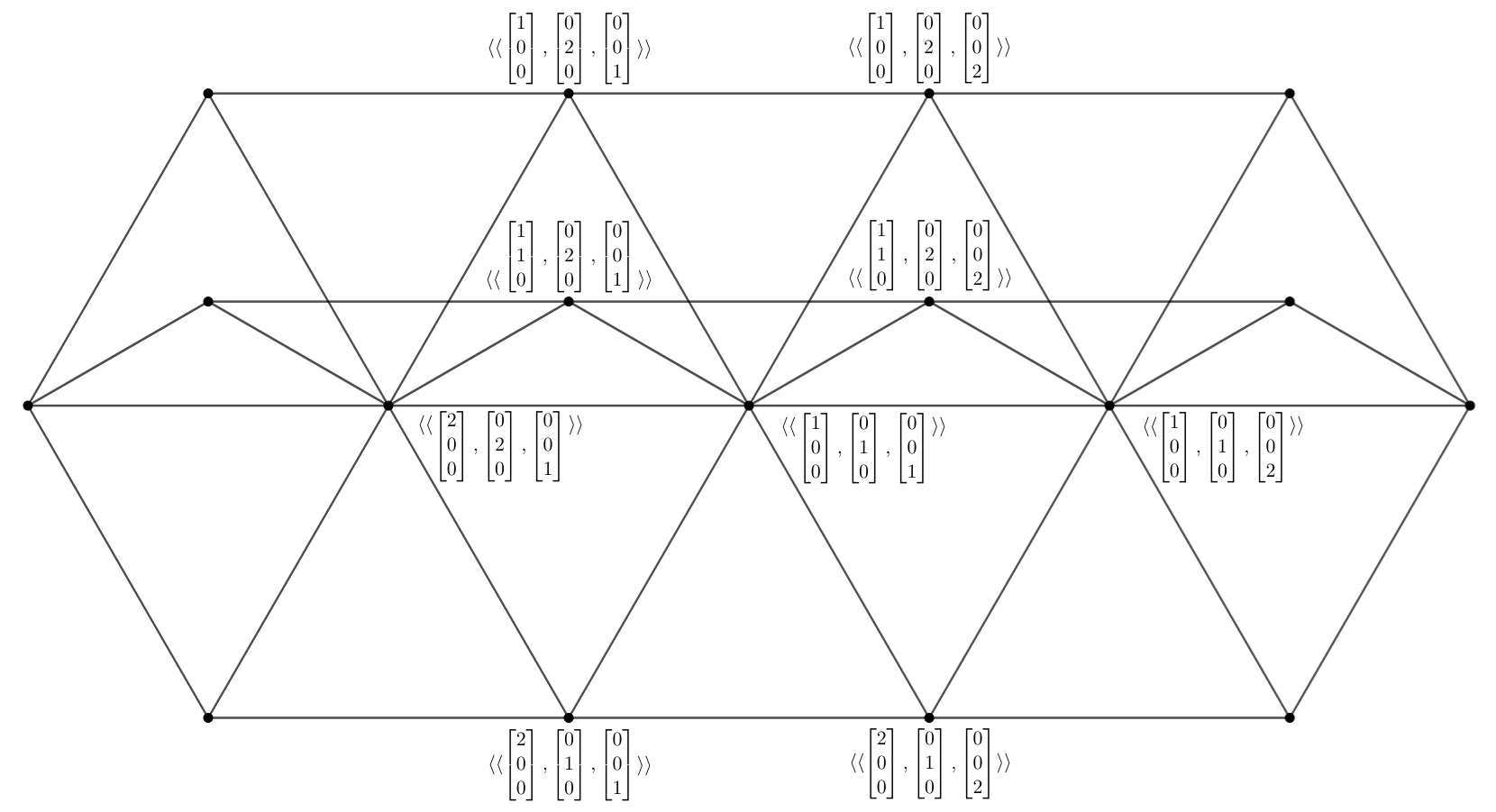}
\caption{Illustration of a part of the building involved in the proof of Proposition \ref{proparithmetic}, when $K=\Q$ and $\nu$ is the $2$-adic valuation on $\Q$. The key feature which makes the argument work is that the linear automorphism $g$ acts as a translation along the illustrated ``slice'' of the building, without permuting the $3$ ``flaps''.}
\label{arithmetic1660}
\end{center}
\end{figure}

\begin{preuve}

Without loss of generality, by transitivity of $\text{Aut}(\widetilde{\text{A}_2}(K,v))$, we may assume that the isomorphism identifies : \[ v_0 = \langle \langle \begin{bmatrix}  1 \\ 0 \\ 0  \end{bmatrix} , \begin{bmatrix}  0 \\ 1 \\ 0  \end{bmatrix},\begin{bmatrix}  0 \\ 0 \\ 1  \end{bmatrix} \rangle \rangle , \] and \[ v_1 =  \langle \langle \begin{bmatrix}  1 \\ 0 \\ 0  \end{bmatrix} ,\begin{bmatrix}  0 \\ 1 \\ 0  \end{bmatrix},\begin{bmatrix}  0 \\ 0 \\ \pi  \end{bmatrix} \rangle \rangle . \]  Then, we note that $\Lambda$ is the set consisting of : \[ L_d = \langle \langle \begin{bmatrix}  1 \\ d \\ 0  \end{bmatrix} ,\begin{bmatrix}  0 \\ \pi \\ 0  \end{bmatrix},\begin{bmatrix}  0 \\ 0 \\ \pi  \end{bmatrix} \rangle \rangle , \] where the coefficients $d\in O_K$ are such that in the reduction $O_K \to O_K/\pi O_K$, the vectors $\begin{bmatrix} 1 , d  \end{bmatrix}$ are pairwise distinct, and : \[   L_\infty = \langle \langle \begin{bmatrix}  \pi \\ 0 \\ 0  \end{bmatrix} ,\begin{bmatrix}  0 \\ 1 \\ 0  \end{bmatrix},\begin{bmatrix}  0 \\ 0 \\ \pi  \end{bmatrix} \rangle \rangle  . \] Let : \[ g = \begin{bmatrix} 1 & 0 & 0 \\ 0 & 1 & 0 \\ 0 & 0 & \pi  \end{bmatrix} \in \text{GL}_3(K) .  \] It induces a linear automorphism of $\widetilde{\text{A}_2}(K,\nu)$ such that $gv_0 = v_1$, and thus it restricts to an isomorphism $\varphi : \Pi_0 \to \Pi_1$. Let $v_2 = \varphi^{-1}(v_0) \in \Pi_0$, and observe that it is opposite to $v_1$. So, for every $h\in H_0$, by Lemma \ref{lemmaextension} we may assume that $h(L_2) = L_2$ without affecting its local action on $\Lambda$. But then : \[ \varphi h \varphi^{-1}(v_0) = \varphi h(v_2) = \varphi(v_2) = v_0, \] hence $\varphi h \varphi^{-1} \in H_1$. Now, the simple computation : \[ \varphi^{-1}(L_d) = g^{-1}L_d = \langle \langle \begin{bmatrix}  1 \\ d \\ 0  \end{bmatrix} ,\begin{bmatrix}  0 \\ \pi \\ 0  \end{bmatrix},\begin{bmatrix}  0 \\ 0 \\ 1  \end{bmatrix} \rangle \rangle ,  \] \[ \varphi^{-1}(L_\infty) = g^{-1}L_\infty = \langle \langle \begin{bmatrix}  \pi \\ 0 \\ 0  \end{bmatrix} ,\begin{bmatrix}  0 \\ 1 \\ 0  \end{bmatrix},\begin{bmatrix}  0 \\ 0 \\ 1  \end{bmatrix} \rangle \rangle  , \] shows that for all $L\in \Lambda$, we have that $\varphi^{-1}(L)$ is the unique $L' \in \Pi_0$ such that it is adjacent to $L$ and $v_2$. But $h$ fixes $v_2$ so $h\varphi^{-1}(L)$ is adjacent to $h(L)$ and $v_2$, and thus $h\varphi^{-1}(L) = \varphi^{-1}h(L)$. In other words, for all $L\in \Lambda$, we have that $\varphi h \varphi^{-1}(L) = h(L)$, and thus the image of $\varphi h \varphi^{-1}$ in $Q_1$ is precisely equal to the image of $h$ in $Q_0$. See Figure \ref{arithmetic1660}. $\square$

\end{preuve}

\begin{prop}
\label{propconj}

Let : \[ \varphi : a_{x}^{\pm 1} \in \Pi_0 \mapsto a_0^{-1} a_{-x}^{\mp 1} \in \Pi_1 . \] Then, $\varphi : \Pi_0 \to \Pi_1$ is an isomorphism such that : \[  \varphi(a_s) = a_{\lambda^{\kappa(o(s))}(s)} , \] for all $s\in S$. In particular, if $\sigma \in \emph{Sym}(\Lambda)$ denotes the restriction of $\varphi$ to $\Lambda$, then $Q_1 = \sigma Q_0 \sigma^{-1}$ and thus $Q_1 = Q_0$ if and only if $\sigma \in N_{\emph{Sym}(\Lambda)}(Q_0)$.

\end{prop}

\begin{figure}
\begin{center}
\includegraphics[scale=0.55]{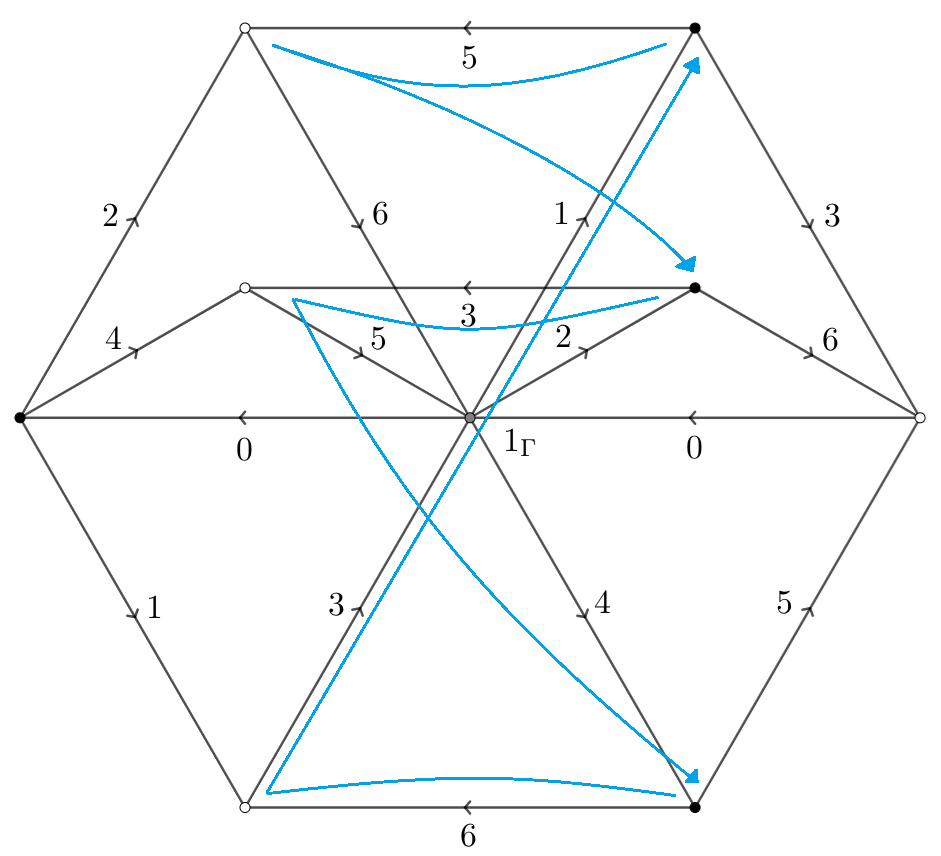}
\caption{Illustration of a part of the link $\Pi_0$ of $1_\Gamma$ when $q=2$. The blue arrows represent the permutation $\sigma\in \text{Sym}(\Lambda)$ which conjugates $Q_1$ and $Q_0$.}
\label{exotic941}
\end{center}
\end{figure}

\begin{preuve}

We note that $\varphi$ is the composition of the automorphism $\mu : a_{x}^{\pm 1} \mapsto a_{-x}^{\mp 1}$ of $\Pi_0$ (see Proposition \ref{propmu}) followed by the restriction of the automorphism $\gamma \mapsto a_0^{-1} \gamma$ of $\Delta$, so it is indeed an isomorphism of the links. Moreover, for all $s\in S$, we have : \[ \varphi(a_s) = a_0^{-1} a_{-s}^{-1} = a_{\lambda^{\kappa(o(s))}(s)} (a_{-s} a_0 a_{\lambda^{\kappa(o(s))}(s)})^{-1} = a_{\lambda^{\kappa(o(s))}(s)} ,  \] where we have used the fact that : \[ (-s,0,\lambda^{\kappa(o(s))}(s)) = (x,x+s,x+s+\lambda^{\kappa(o(s))}(s)) \in T_\kappa ,  \] for $x=-s \in G$. So $\varphi$ preserves $\Lambda$. Moreover : \[ \varphi(v_1) = \varphi(a_0^{-1}) = a_0^{-1} a_0 = 1_\Gamma = v_0 , \] and thus $H_1 = \varphi H_0 \varphi^{-1}$. Restricting to $\Lambda$, we get $Q_1 = \sigma Q_0 \sigma^{-1}$. See Figure \ref{exotic941} for an illustration. $\square$

\end{preuve}

\begin{lemme}
\label{lemmanorm}

Let $Q = \emph{P$\Gamma$L}_2(\F_q)$ be seen as a subgroup of $\emph{Sym}(q+1)$ via its action on the projective line. Then, $N_{\emph{Sym}(q+1)}(Q) = Q$. In particular, the normaliser has order $e(q-1)q(q+1)$.

\end{lemme}

\begin{preuve}

First, we show that the centraliser of $Q$ is trivial. Note that : \[ U = \set{ \begin{bmatrix} 1 & x \\ 0 & 1  \end{bmatrix} : x \in \mathbb{F}_q } ,  \] appears as a subgroup of $Q$ and acts on the projective line by fixing one point, and regularly on the $q$ remaining ones. In particular, a permutation that centralises $Q$ must centralise $U$ and hence fix the same point. But $Q$ acts transitively on the projective line, and thus a permutation that centralises $Q$ must centralise every conjugate of $U$ and hence it fixes every point on the projective line. This implies that the natural homomorphism $N_{\text{Sym}(q+1)}(Q) \to \text{Aut}(Q)$ is injective. We conclude using the fact that $\text{Aut}(\text{P$\Gamma$L}_2(\F_q))\cong\text{P$\Gamma$L}_2(\F_q)$. $\square$

\end{preuve}

\begin{preuve} 

(of Theorem \ref{theomain}) If $\Delta$ is isomorphic to a Bruhat-Tits building, then by Proposition \ref{proparithmetic} we must have $Q_1 = Q_0$. By Proposition \ref{propconj}, we have $Q_1 = \sigma Q_0 \sigma^{-1}$ ; and by assumption, $\sigma \notin Q_0$, and thus by Lemma \ref{lemmanorm}, we have that $\sigma$ does not normalise $Q_0$. Thus $Q_1\neq Q_0$, a contradiction. The second part of the statement on the number of quasi-isometry classes then follows from the fact that $q^2+q+1 = |G| \leq |\text{Aut}(T_\kappa)|$ and $|\text{Aut}(F(G,S))| \leq |\text{Aut}(\text{A}_2(\F_q))| = e(q^2+q+1)(q-1)^2q^3(q+1)$, and applying Proposition \ref{propcount} and Proposition \ref{propqi}. $\square$

\end{preuve}


\section{More examples and other link types}
\label{sec4}

Let $q\geq 2$ be a prime power.

\begin{prop}
\label{propquad}

Let $G = \Z/(q^4+q^2+1)$ and let $S\subseteq G$ be as in Section \ref{sec3}, in such a way that $F(G,S)$ is a symmetric representation of $\emph{A}_2(\F_{q^2})$. Then, there are at least $2^{(q^2-q+1)R(q)}$ distinct triangle presentations compatible with $F(G,S)$.

\end{prop}

\begin{preuve}

Note that : \[ (q^4+q^2+1) = (q^2+q+1)(q^2-q+1) , \] therefore let $H\leq G$ be the subgroup of order $q^2+q+1$. In the isomorphism $G\cong \F_{q^6}^\times/\F_{q^2}^\times$, it corresponds to the image of the subgroup $\F_{q^3}^\times \leq \F_{q^6}^\times$. If $x\in \F_{q^3}^\times$ is such that : \[ \text{Tr}_{\F_q}(x) = 0, \] then seen as an element of $\F_{q^6}$ we can compute : \[ \text{Tr}_{\F_{q^2}}(x) = x+x^{q^2}+x^{q^4} = x+x^q+x^{q^2} = \text{Tr}_{\F_q}(x) = 0, \] where we have used the fact that $x^{q^3}=x$. This shows that $S\cap H$ contains $q+1$ elements and is preserved by $\lambda : x\in G \mapsto q^2 x \in G$, because it restricts to the (square of the) Frobenius automorphism of $\F_{q^3}$ over $\F_q$. The counting of Proposition \ref{proporbit} shows that $S\cap H$ contains $R(q)$ orbits of length $3$, and we may apply Proposition \ref{propbreak} to conclude. $\square$

\end{preuve}

\begin{ex}

For $q=2$, Proposition \ref{propquad} refers to the triangle presentations which we have already described in Example \ref{exquad}.

\end{ex}

The graph which we denote by $\text{Opp}(\text{A}_2(\F_q))$ is the subgraph of $\text{A}_2(\F_q)$ generated by the vertices which are opposite to : \[  V_1 = \langle \begin{bmatrix} 1 \\ 0 \\ 0 \end{bmatrix} \rangle  \text{ or } V_2 = \langle \begin{bmatrix} 1 \\ 0 \\ 0   \end{bmatrix} , \begin{bmatrix} 0 \\ 1 \\ 0   \end{bmatrix} \rangle . \] It is isomorphic to the coset graph of : \[ U = \set{ \begin{bmatrix} 1 & x & z \\ 0 & 1 & y \\ 0 & 0 & 1  \end{bmatrix} : x,y,z\in \F_q }  , \] with respect to the subgroups : \[ U_1 = \set{ \begin{bmatrix} 1 & x & 0 \\ 0 & 1 & 0 \\ 0 & 0 & 1  \end{bmatrix} : x\in \F_q } ,\] \[ U_2  = \set{ \begin{bmatrix} 1 & 0 & 0 \\ 0 & 1 & y \\ 0 & 0 & 1  \end{bmatrix} : y \in \F_q } . \] Combining these two points of view, it is possible to deduce a number of informations on $\text{Opp}(\text{A}_2(\F_q))$, namely : \begin{itemize} \item It has $2q^2$ vertices and $q^3$ edges. \item It is regular of degree $q$. \item It is bipartite. \item It is connected. \item It has girth $6$ if $q\geq 3$ and $8$ if $q=2$. \item It has diameter $4$. \item It has spectral gap $1-\frac{\sqrt{q}}{q}$. \end{itemize} Let further : \[  G = \set{ \begin{bmatrix} 1 & y & z \\ 0 & 1 & y \\ 0 & 0 & 1  \end{bmatrix} : y,z \in \F_q }  , \] and let : \[ S = \set{ \begin{bmatrix} 1 & y & y^2 \\ 0 & 1 & y \\ 0 & 0 & 1  \end{bmatrix} : y \in \F_q } . \] The group $G$ is isomorphic to the additive group $(\F_q)^2$, but this matrix representation will be useful in proving the following proposition.

\begin{prop}

The graph $\mathcal{G}^{(q^2)}_{F(G,S)}$ is isomorphic to $\emph{Opp}(\emph{A}_2(\F_q))$.

\end{prop}

\begin{preuve}

We start by providing a numbering of the vertices of $\text{Opp}(\text{A}_2(\F_q))$, which we think of cosets in $U/U_1 \cup U/U_2$. Note that $U/U_1 = \set{ g_{y,z}U_1 : y,z\in \F_q}$, where we use the notation : \[ g_{y,z} = \begin{bmatrix} 1 & y & z \\ 0 & 1 & y \\ 0 & 0 & 1  \end{bmatrix}  \in G.  \] Indeed, if $g\in U$ is an arbitrary element of the form : \[ g = \begin{bmatrix} 1 & x & z \\ 0 & 1 & y \\ 0 & 0 & 1  \end{bmatrix} , \] then $gU_1 = g_{y,z}U_1$. Similarly for $U/U_2$, we have that $U/U_2 = \set{g_{y,z}U_2 : y,z \in \F_q}$ because $gU_2 = g_{x,z+(x-y)x}U_2$. With this numbering, we see that : \[ g_{y_1,z_1}U_1 \cap g_{y_2,z_2} U_2 \neq \emptyset \text{ if and only if } -z_1+z_2-y_2(-y_1+y_2) = 0 \text{ if and only if } g_{y_1,z_1}^{-1}g_{y_2,z_2} \in S ,   \] which ends our proof. $\square$

\end{preuve}

\begin{figure}[h!]
\begin{center}
\includegraphics[scale=0.4]{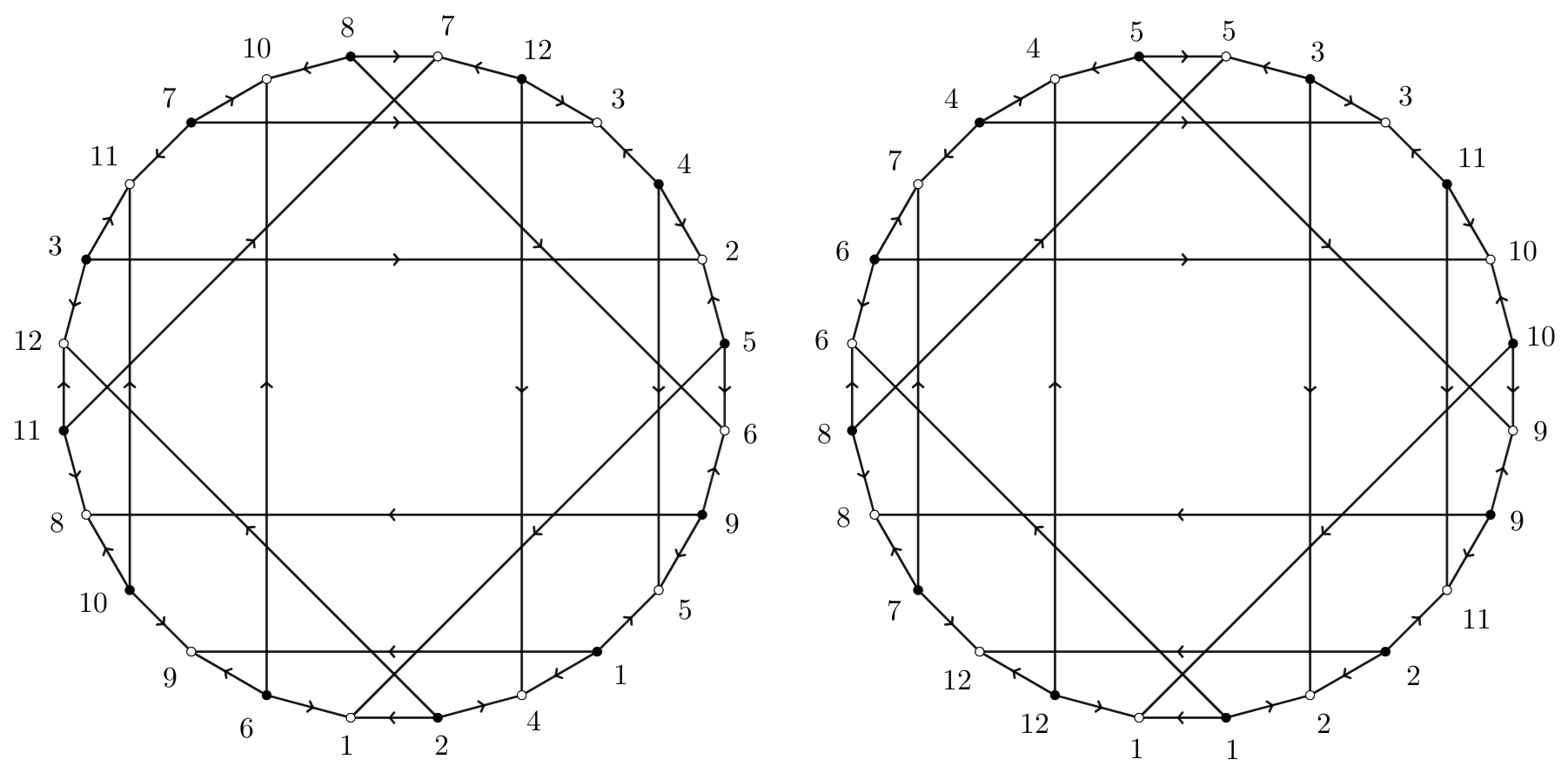}
\caption{Two numberings of the Nauru graph, which corresponds to subsets $F_1,F_2\subseteq \set{1,\dots,12}^2$ on the left and on the right, respectively.}
\label{nauru1802}
\end{center}
\end{figure}

\begin{theo}
\label{theoopp}

If $q=1 \emph{ mod } 3$, let $\alpha \in \F_q$ be an element of order $3$. Then, the map defined by : \[ \lambda(\begin{bmatrix} 1 & y & y^2 \\ 0 & 1 & y \\ 0 & 0 & 1 \end{bmatrix}) = \begin{bmatrix} 1 & \alpha y & (\alpha y)^2 \\ 0 & 1 & \alpha y \\ 0 & 0 & 1 \end{bmatrix}  ,  \] for all $y\in \F_q$, defines a permutation $\lambda : S \to S$ of order $3$ such that : \[ s\lambda(s)\lambda^2(s) = s\lambda^2(s)\lambda(s) = \begin{bmatrix} 1 & 0 & 0 \\ 0 & 1 & 0 \\ 0 & 0 & 1 \end{bmatrix} , \] for all $s\in S$. In particular, there exist at least $2^{\frac{q-1}{3}}$ distinct triangle presentations compatible with $F(G,S)$.

\end{theo}

\begin{preuve}

It is a direct matrix computation and application of Proposition \ref{propbreak} with $H=G$. $\square$

\end{preuve}

\begin{ex}
\label{exnauru}

We have illustrated two representations of the Nauru graph of the form $\mathcal{G}^{(12)}_F$ in Figure \ref{nauru1802}. It has girth $6$ and diameter $4$. It supports the triangle presentations $T_1$ and $T_2$ given in Tables $4$ and $5$, and we have found using Magma that $\Gamma^{(12)}_{T_1}$ is a hyperbolic group while $\Gamma^{(12)}_{T_2}$ contains a subgroup isomorphic to $\Z^2$ (for instance, the subgroup generated by $a_5a_1$ and $a_{10}a_4$). We see this as evidence that some of the groups that arise from the triangle presentations of Theorem \ref{theoopp} could be hyperbolic, because the link also has girth $6$ and diameter $4$ in this case.

\end{ex}


\section{Tables}

\begin{table}[h!]
\label{tablequad1:1}
\caption{The triangle presentation $T_1$ from Example \ref{exquad}.}
\begin{center}
\begin{tabular}{ |p{2cm}|p{2cm}|p{2cm}|p{2cm}|p{2cm}| } 
\hline
(0,7,14) & (0,9,3) & (0,14,7) & (0,15,12) & (0,18,6) \\
(1,8,15) & (1,10,4) & (1,15,8) & (1,16,13) & (1,19,7) \\
(2,9,16) & (2,11,5) & (2,16,9) & (2,17,14) & (2,20,8) \\
(3,0,9) & (3,10,17) & (3,12,6) & (3,17,10) & (3,18,15) \\
(4,1,10) & (4,11,18) & (4,13,7) & (4,18,11) & (4,19,16) \\
(5,2,11) & (5,12,19) & (5,14,8) & (5,19,12) & (5,20,17) \\
(6,0,18) & (6,3,12) & (6,13,20) & (6,15,9) & (6,20,13) \\
(7,0,14) & (7,1,19) & (7,4,13) & (7,14,0) & (7,16,10) \\
(8,1,15) & (8,2,20) & (8,5,14) & (8,15,1) & (8,17,11) \\
(9,2,16) & (9,3,0) & (9,6,15) & (9,16,2) & (9,18,12) \\
(10,3,17) & (10,4,1) & (10,7,16) & (10,17,3) & (10,19,13) \\
(11,4,18) & (11,5,2) & (11,8,17) & (11,18,4) & (11,20,14) \\
(12,0,15) & (12,5,19) & (12,6,3) & (12,9,18) & (12,19,5) \\
(13,1,16) & (13,6,20) & (13,7,4) & (13,10,19) & (13,20,6) \\
(14,0,7) & (14,2,17) & (14,7,0) & (14,8,5) & (14,11,20) \\
(15,1,8) & (15,3,18) & (15,8,1) & (15,9,6) & (15,12,0) \\
(16,2,9) & (16,4,19) & (16,9,2) & (16,10,7) & (16,13,1) \\
(17,3,10) & (17,5,20) & (17,10,3) & (17,11,8) & (17,14,2) \\
(18,4,11) & (18,6,0) & (18,11,4) & (18,12,9) & (18,15,3) \\
(19,5,12) & (19,7,1) & (19,12,5) & (19,13,10) & (19,16,4) \\
(20,6,13) & (20,8,2) & (20,13,6) & (20,14,11) & (20,17,5) \\
\hline
\end{tabular}
\end{center}
\end{table}

\begin{table}
\label{tablequad2:1}
\caption{The triangle presentation $T_2$ from Example \ref{exquad}. In red, we highlight the triples that differ from $T_1$.}
\begin{center}
\begin{tabular}{ |p{2cm}|p{2cm}|p{2cm}|p{2cm}|p{2cm}| } 
\hline
(0,7,14) & (0,9,3) & (0,14,7) & (0,15,12) & (0,18,6) \\
(1,8,15) & (1,10,4) & (1,15,8) & (1,16,13) & (1,19,7) \\
(2,9,16) & \textcolor{red}{(2,11,8)} & (2,16,9) & \textcolor{red}{(2,17,5)} & \textcolor{red}{(2,20,14)} \\
(3,0,9) & (3,10,17) & (3,12,6) & (3,17,10) & (3,18,15) \\
(4,1,10) & (4,11,18) & (4,13,7) & (4,18,11) & (4,19,16) \\
\textcolor{red}{(5,2,17)} & (5,12,19) & \textcolor{red}{(5,14,11)} & (5,19,12) & \textcolor{red}{(5,20,8)} \\
(6,0,18) & (6,3,12) & (6,13,20) & (6,15,9) & (6,20,13) \\
(7,0,14) & (7,1,19) & (7,4,13) & (7,14,0) & (7,16,10) \\
(8,1,15) & \textcolor{red}{(8,2,11)} & \textcolor{red}{(8,5,20)} & (8,15,1) & \textcolor{red}{(8,17,14)} \\
(9,2,16) & (9,3,0) & (9,6,15) & (9,16,2) & (9,18,12) \\
(10,3,17) & (10,4,1) & (10,7,16) & (10,17,3) & (10,19,13) \\
(11,4,18) & \textcolor{red}{(11,5,14)} & \textcolor{red}{(11,8,2)} & (11,18,4) & \textcolor{red}{(11,20,17)} \\
(12,0,15) & (12,5,19) & (12,6,3) & (12,9,18) & (12,19,5) \\
(13,1,16) & (13,6,20) & (13,7,4) & (13,10,19) & (13,20,6) \\
(14,0,7) & \textcolor{red}{(14,2,20)} & (14,7,0) & \textcolor{red}{(14,8,17)} & \textcolor{red}{(14,11,5)} \\
(15,1,8) & (15,3,18) & (15,8,1) & (15,9,6) & (15,12,0) \\
(16,2,9) & (16,4,19) & (16,9,2) & (16,10,7) & (16,13,1) \\
(17,3,10) & \textcolor{red}{(17,5,2)} & (17,10,3) & \textcolor{red}{(17,11,20)} & \textcolor{red}{(17,14,8)} \\
(18,4,11) & (18,6,0) & (18,11,4) & (18,12,9) & (18,15,3) \\
(19,5,12) & (19,7,1) & (19,12,5) & (19,13,10) & (19,16,4) \\
(20,6,13) & \textcolor{red}{(20,8,5)} & (20,13,6) & \textcolor{red}{(20,14,2)} & \textcolor{red}{(20,17,11)} \\
\hline
\end{tabular}
\end{center}
\end{table}

\begin{table}
\label{tableheawood}
\caption{The triangle presentation $T$ from Example \ref{exheawood}.}
\begin{center}
\begin{tabular}{ |p{2cm}|p{2cm}|p{2cm}|} 
\hline
(0,1,3) & (0,2,6) & (0,4,5) \\
(1,2,4) & (1,3,0) & (1,5,6) \\
(2,3,5) & (2,4,1) & (2,6,0)  \\
(3,0,1) & (3,4,6) & (3,5,2) \\
(4,1,2) & (4,5,0) & (4,6,3) \\
(5,0,4) & (5,2,3) & (5,6,1) \\
(6,0,2) & (6,1,5) & (6,3,4) \\
\hline
\end{tabular}
\end{center}
\end{table}

\begin{table}
\label{tablenauru1}
\caption{The triangle presentation $T_1$ from Example \ref{exnauru}.}
\begin{center}
\begin{tabular}{|p{2cm}|p{2cm}|p{2cm}|} 
\hline
(1,2,4) & (1,5,9) & (1,6,5) \\
(2,3,12) & (2,4,1) & (2,5,4) \\
(3,4,12) & (3,7,11) & (3,12,2) \\
(4,1,2) & (4,2,5) & (4,12,3) \\
(5,1,6) & (5,4,2) & (5,9,1) \\
(6,5,1) & (6,8,9) & (6,9,10) \\
(7,8,10) & (7,11,3) & (7,12,11) \\
(8,9,6) & (8,10,7) & (8,11,10) \\ 
(9,1,5) & (9,6,8) & (9,10,6) \\
(10,6,9) & (10,7,8) & (10,8,11) \\
(11,3,7) & (11,7,12) & (11,10,8) \\
(12,2,3) & (12,3,4) & (12,11,7) \\
\hline
\end{tabular}
\end{center}
\end{table}

\begin{table}
\label{tablenauru2}
\caption{The triangle presentation $T_2$ from Example \ref{exnauru}.}
\begin{center}
\begin{tabular}{ |p{2cm}|p{2cm}|p{2cm}| } 
\hline
(1,1,1) & (1,10,6) & (1,12,2) \\
(2,1,12) & (2,2,2) & (2,3,11) \\ 
(3,3,3) & (3,4,5) & (3,11,2) \\
(4,4,4) & (4,5,3) & (4,12,7) \\
(5,3,4) & (5,5,5) & (5,8,9) \\
(6,1,10) & (6,6,6) & (6,8,7) \\
(7,4,12) & (7,6,8) & (7,7,7) \\
(8,7,6) & (8,8,8) & (8,9,5) \\
(9,5,8) & (9,9,9) & (9,10,11) \\
(10,6,1) & (10,10,10) & (10,11,9) \\
(11,2,3) & (11,9,10) & (11,11,11) \\
(12,2,1) & (12,7,4) & (12,12,12) \\
\hline
\end{tabular}
\end{center}
\end{table}


\newpage

\bibliography{ref}{}
\bibliographystyle{plain}


\end{document}